\newif\ifdraft
\draftfalse

\ifdraft
\documentclass{amsart}
\else
\documentclass[12pt]{amsart}
  \usepackage{a4wide}
 \fi
\usepackage{mathrsfs}
\usepackage{amssymb}                    
\usepackage{epic}
\usepackage{eepic}

 \date{2004-04-11}

\newcommand{\details}[1]{}

\newcommand{\power}{{\mathscr P}}
\newcommand{\X}{{\mathscr X}}
\newcommand{\WF}{{\bf WF}}
\newdimen\dick
\newdimen\hoch
\newdimen\breit
\newdimen\hhoch
\newdimen\vhoch

\hhoch=-\hoch
\advance\hhoch\dick

\vhoch=\hoch
\advance\vhoch-0.5\dick

\newcommand{\T}{{\sf T}}

\newcommand{\C}{\mathscr{C}}

\newcommand{\F}{\mathscr{F}}
\newcommand{\M}{{\mathscr M}}
\newcommand{\N}{{\mathbb N}}

\newcommand{\ran}{{\rm ran}}
\newcommand{\fullclone}{{\mathscr O}}

\makeatletter
\newcommand{\KNUTHcases}[1]
{\left \{\,\vcenter {\normalbaselines \m@th
\ialign {$##\hfil $&\quad ##\hfil \crcr #1\crcr }}\right .
}  
\makeatother

\def \itm#1 {\item[(#1)]}

\newcommand{\on}{{\upharpoonright}}
\newcommand{\oo}[1]{\fullclone^{(#1)}}

\newcommand{\nd}{\mathord{\raise1pt\hbox{$\nabla$}\!\!\Delta}}

\newcommand{\Pol}{{\rm Pol}}

\swapnumbers
\theoremstyle{remark}

\def\xx#1 {\newtheorem{#1}[thm]{#1}}
\xx Acknowledgement
\xx Lemma
\xx {Main Lemma}
\xx Theorem
\xx Corollary
\xx Idea
\xx {First Step}
\xx Consequence
\xx Example
\xx Definition
\xx Fact
\xx {Definition and Fact}
\xx{Open Questions}
\xx Setup
\xx Notation
\xx {More Notation}
\xx Remark 
\xx Conclusion
\xx Question
\title{Analytic Clones}

\def\pre#1^#2 {{}^{#2} #1}

\newcommand{\dom}{{\rm dom}}


\author{Martin Goldstern}

\address{Discrete Mathematics and Geometry\\
 Algebra Research Group\endgraf TU Wien
\\
Wiedner Hauptstra\3e 8--10/104
\endgraf
1040 Wien, Austria (Europe)}
\email{Martin.Goldstern@tuwien.ac.at}
\urladdr{http://www.tuwien.ac.at/goldstern/}

\begin{document}

\begin{abstract}
We use a method from descriptive set theory to investigate 
the two complete clones above the unary clone on a countable set.
\end{abstract}

\maketitle

\section{Introduction.  Known results}

An ``operation'' on a set $X$ is a function $f:X^n\to X$, for some $n\in \N\setminus\{0\}$. 
If $f$ is such an  ``$n$-ary operation'', $g_1$, \ldots, $g_n$ are
$k$-ary, then the ``composition'' $f(g_1,\ldots, g_n)$ is defined 
naturally:
$$
f(g_1,\ldots, g_n)(\vec x) = 
f(g_1(\vec x), \ldots, g_n(\vec x)  ) \ \ \mbox{ for all $\vec x \in X^k$}
$$ 
 A clone on a set $X$ is a set $\C$ of  operations
 which contains all the projections and is closed under 
 composition.
 (Alternatively, $\C$ is a clone on~$X$
 {if} $\C$ is the set of term functions of 
 some universal  algebra over~$X$.)  
  
  The family of all clones forms a complete algebraic 
 lattice $Cl(X)$ with greatest element $\fullclone = \bigcup_{n=1}^\infty \oo n $,  
 where $\oo n = X^{X^n} $ is the set of all $n$-ary operations on~$X$.
   (In this paper, the underlying set $X$  will always be the set 
 $\N = \{0,1,2,\ldots \} $ of natural numbers.) 
  
   The coatoms of this lattice $Cl(X) $ are called 
 ``precomplete clones'' or ``maximal clones'' on~$X$.  
  
For any set $\C \subseteq \fullclone$ we write $\langle \C\rangle$ 
 for the
smallest
clone containing $\C$.  In particular, $\langle \oo 1 \rangle$ is
 the set of all functions $\pi\circ f$, where $f:X\to X$ is arbitrary
 and 
$\pi:X^n\to X$ is a projection to one coordinate.  However, to lighten the
 notation we will identify $\oo 1 $ (the set of all unary functions)
with $\langle \oo 1 \rangle$ (the set of all ``essentially'' unary
 functions).

 For singleton sets $X$ the lattice $Cl(X)$ is trivial; for $|X|=2$
 the lattice $Cl(X)$ is countable, and well understood (``Post's
 lattice'').  For $|X|\ge 3$, $Cl(X)$ is uncountable. 
 For infinite $X$, $Cl(X)$ has $2^{2^{|X|}}$ elements, and there are
 even $2^{2^{|X|}}$ precomplete clones on $X$.

In this paper we are interested in the interval $[\oo 1, \fullclone]$ of the
clone lattice on a countable set $X$.  It will turn out that 
methods from descriptive set theory are useful to describe the 
complexity of several interesting  clones in this interval, and also
the overall structure of the interval.

For simplicity we concentrate on binary clones, i.e., clones generated
by binary functions.  Equivalently, we can define a binary clone to be
a set $\C$ of functions $f:\N^2\to \N$ which contains the two
projections and is closed under composition: if $f,g,h\in \C$, then
also the function $f(g,h)$ (mapping  
 $(x,y)$ to $ f(g(x,y),h(x,y))$) is in $\C$. 

The set of binary clones, $Cl^{(2)}(X)$, also forms a complete
 algebraic lattice.

Occasionally we will remark on how to modify the definitions or
theorems for the case of ``full''  clones, i.e., for clones that are
not necessarily generated by binary functions.  (In some cases this
generalization is trivial, in other cases it is nontrivial but known,
and in some cases it is still open.)

By \cite{Gavrilov:1965} (see also \cite{GoSh:737}), we know that there
are exactly 2 precomplete binary clones above $\oo 1 $, which we call
$\T_1$ and $\T_2$ (see below).  
 It is known that the interval 
$[\oo 1, \oo 2 ] $ of binary clones 
is dually atomic, so it can be written 
$$ [\oo 1, \oo 2 ]  = [\oo 1 , \T_1] \cup [\oo 1 , \T_2]\cup \{\oo 2 \},  $$
i.e., every binary clone above $\oo 1 $ other than $\oo 2 $ itself 
is contained in $\T_1$ or in $\T_2$.

  So we will have to investigate the
intervals $ [\oo 1, \T_1]$ and $[\oo 1, \T_2]$.  We will see that these
two structures are very different, and that this difference can be
traced back to a difference in ``complexity'' of the two binary clones
$\T_1$ and $\T_2$.  

More precisely, $\T_1$ is a Borel set, while $\T_2$ is a complete
coanalytic set.  We will see that $\T_1$ is finitely
generated over $\oo 1$, but $\T_2$ is not countably generated over 
$\oo 1 $.

\begin{Definition} 
A function $f:\N\times \N$ is called ``almost unary'', if at least one
of the following holds:
\begin{enumerate}
\item[(x)] There is a function $F:\N\to \N$ such that $\forall
             x\,\forall y: f(x,y)\le F(x)$.  
\item[(y)] There is a function $F:\N\to \N$ such that $\forall
             x\,\forall y: f(x,y)\le F(y)$.  
\end{enumerate}
We let $\T_1$ be the set of all binary functions which are almost
unary. It is easy to see that $\T_1$ is a binary clone containing
$\oo 1 $. 
\end{Definition}

\begin{Definition}
Let $B \subseteq \oo 2$. The set $\Pol(B)$ is defined as 
$$
\bigcup_{k=1 }^\infty \{f\in \oo k : \forall  g_1,\ldots, g_k\in B:
              f(g_1,\ldots, g_k)\in B\}$$

(Background and a more general definition of $\Pol$ can be found in 
\cite{PK:1979}.)

\end{Definition}
\begin{Fact}
$\Pol(B)$ is a clone.  If $B$ is a binary clone, then $\Pol(B) 
\cap \oo 2 = B$. 
\end{Fact}

\begin{Definition}
Let $\Delta:= \{(x,y)\in \N\times \N: x > y\}$, $\nabla:= \{(x,y): x < y\}$. 

For  $S_1,S_2\subseteq \N$ we let $\Delta_{S_1,S_2}: =\Delta\cap
(S_1\times S_2)$.  We define $\nabla_{S_1,S_2}$
similarly.

If $S_1,S_2$ are infinite subsets of $\N$, and $g:
\Delta_{S_1,S_2}\to \N$ or $g:\nabla_{S_1,S_2} \to \N$, then we say
that $g$ is ``canonical'' iff one of the following holds: 
\begin{enumerate}
\item $g$ is constant
\item There is a 1-1 function $G:S_1\to \N$ such that 
$$ \forall (x,y)\in \dom(g): g(x,y)  = G(x) $$
\item There is a 1-1 function $G:S_2\to \N$ such that 
$$ \forall (x,y)\in \dom(g): g(x,y)  = G(y) $$
\item $g$ is 1-1. 
\end{enumerate}
The ``type'' of $g$ is one  of the labels
``constant'', ``x'', ``y'', or ``1-1'', respectively. 

Let $f:\N\times \N \to \N$.    We say that $f$ is 
canonical on $S_1\times S_2$ 
iff both functions $f\on  \Delta_{S_1,S_2}$ and
$f\on \nabla_{S_1,S_2}$ are canonical (but not necessarily
of the same type), and moreover: 
\begin{quote} 
Either  the ranges of 
 $f\on  \Delta_{S_1,S_2}$ and $f\on \nabla_{S_1,S_2}$ are disjoint,
\\
or $S_1=S_2$, and $f(x,y)=f(y,x)$ for all $x,y\in S_1$. 
\end{quote}


\end{Definition}

The following fact is a consequence of Ramsey's theorem, see 
\cite{GoSh:737}.   It was originally proved in a slightly different formulation
already in   \cite{Gavrilov:1965}.

\begin{Fact}\label{ramsey}
 Let $f:\N\times \N\to \N$. Then there are   infinite sets
$S_1$, $S_2$ such that $f$ is canonical on $S_1\times S_2$. 

Moreover, for any infinite sets $S_1, S_2$ we can find 
infinite $S'_1\subseteq S_1 $, $S'_2\subseteq S_2$
 such that $f$ is canonical on $S'_1\times S'_2$. 
\end{Fact}

\begin{Definition}
Let $f:\N\times \N \to \N$.  We say that 
$f$ is              
 ``nowhere injective'', if: 
\begin{quote}  whenever $f$ is canonical on $S_1\times S_2$, then
neither $f\on \Delta_{S_1,S_2}$ nor $f\on \nabla_{S_1,S_2}$  is 1-1. 
\end{quote}

We let $\T_2$ be the set of all nowhere injective functions.   Using
fact~\ref{ramsey}, it is easy to check 
 that $\T_2$  is a binary clone; clearly $\T_2$  contains $\oo 1$. 
(More precisely, $\T_2$ contains $\langle \oo 1 \rangle\cap \oo 2$.)

\end{Definition}

\begin{Theorem}[Gavrilov \cite{Gavrilov:1965}]
$\T_1$ and $\T_2$ are precomplete binary clones, and every binary clone
containing $\oo 1$ is either contained in one of $\T_1$, $\T_2$, or
equal to the clone of all binary functions. 

[For the non-binary case: $\Pol(\T_1)$ and $\Pol(\T_2)$ are precomplete
clones, and every clone $\supseteq \oo 1$ is either $=\fullclone$, or
$\subseteq \Pol(\T_1)$, or $\subseteq \Pol(\T_2)$.]
\end{Theorem}

We will prove the following:

\begin{itemize}
\item (Section \ref{t1}) $\T_1$ is finitely generated over $\oo 1 $, so the
interval $[\oo 1, \T_1]$ in the lattice of binary clones is dually atomic.
\\
In fact, the interval contains a unique coatom: $\T_1\cap \T_2$. 
\item (Section \ref{t2}) 
$\T_2$ is neither finitely nor countably generated
over  $\oo 1 $. \\
$\T_1\cap \T_2$ is a coatom in the interval $[\oo 1, \T_2]$. (Easy)
\\
Any clone which is a Borel set (or even an analytic set) cannot be 
a coatom in this interval. 
\end{itemize}

\begin{Acknowledgement}
I am grateful to 
        J.~Jezek, 
        K.~Kearnes,
        R.~P\"oschel,
        A.~Romanowska,
        A.~Szendrei and R.~Willard for inviting me to 
Bela Csakany's birthday conference (Szeged, 2002), at which I first
presented the main ideas from this paper.
\end{Acknowledgement}

\section{Descriptive Set Theory}
\label{dst}

We collect a few facts and notions from descriptive set theory.
(For motivation, history, details and proofs see the 
textbooks by Moschovakis \cite{Moschovakis:1980} or 
Kechris \cite{Kechris:1995}.)

  Let
$X$ be a countable set (usually $X=\N$, or $X=\N^k$),
and  $Y$ a finite or  countable set, $|Y|\ge
2$ (usually $Y=\N$, or $Y= 2:= \{0,1\}$).
  $Y^X$ is the space of all functions
from $X$ to~$Y$.  

We equip $Y$ with the discrete topology,  $Y^X$ and $Y^{X^n}$ with the product
topology, and $\bigcup_{n=1}^\infty  Y^{X^n}$ with the sum topology. All 
these spaces are ``Polish spaces'', i.e., they are
 separable and carry a (natural) complete metric.

The family of Borel sets is the smallest family $B$
      that contains all open sets and is closed under  complements and
      countable unions (equivalently: contains all open sets and 
all closed sets, and is  closed under countable unions
      and countable intersections).

A  function $f$
 between two topological spaces is called a Borel function iff the preimage
of any Borel set under $f$ is again a Borel set.

A {\em finite sequence on $Y$} is a tuple  $(a_0,\ldots, a_{n-1})\in
Y^n$.   If $s\in Y^k$ and $t\in Y^n$ are finite sequences, $k<n$, then
 we write $s\vartriangleleft t$ iff $s$ is an initial segment
 of~$t$.

We write $Y^{<\omega} := \bigcup_{n\in \N} Y^n$ for the set of all
finite sequences on~$Y$.   If $Y$ is countable, then also 
$Y^{<\omega}$ is countable.  

We can identify $\power(Y^{<\omega})$,
the power set of $Y^{<\omega}$, with the set $2^{Y^{<\omega}}$ of
all characteristic functions, so also 
 $\power(Y^{<\omega})$ carries a natural topology.

A ``tree on $Y$''
 is a set $T  \subseteq \bigcup_{n\in \N}  Y^n$ of finite sequences 
which is  downward closed, i.e., 
\begin{quote} whenever $t\in T$, $s\vartriangleleft t$, 
then also $s\in T$
\end{quote}

The set of all trees is easily seen to be a closed subset
of $\power(Y^{<\omega})$.

For any tree $T$ on~$Y$ we call $f\in Y^\N$ a {\em branch} of~$T$ iff
$\forall n: f\on n \in T$.  (Here we write 
$f\on n$ for $(f(0), \ldots, f(n-1))$.)

We write $[T]$ for the set of all branches of~$T$. 

It is easy to see that $[T]$ is always a closed set in $Y^\N$, and that 
every closed set $\subseteq Y^\N$ is of the form $[T]$ for some tree $T$.

We call a tree $Y$ {\em well-founded} if $[T] = \emptyset$, 
i.e., if there is no sequence $s_0 \vartriangleleft s_1 \vartriangleleft
\cdots$
of elements of~$T$. 

We write ${\WF}$ for the set of all well-founded trees.

The class of {\em analytic sets} is a proper extension of the class of
Borel sets.  There are several possible equivalent definitions of 
``analytic'', for example one could  choose the equivalence 
(1)$\Leftrightarrow$(3) in fact~\ref{def.ana} as the definition of
``analytic''. 

\begin{Fact} \label{def.ana}
 Let $\X$ be a Polish (=complete metric separable) topological 
space, $A \subseteq \X$, $C:= \X \setminus A$. 
Then the following are equivalent: 
\begin{enumerate}
\item $A$ is analytic
\item $C$ is coanalytic
\item $A=\emptyset$, or there is a continuous function $f:\N^\N\to \X$
with $A=f[\N^\N]$
\item There is a Borel set $B \subseteq \N^\N$ and a continuous function 
 $f:\N^\N\to\X$ with $A=f[B]$
\item There is a continuous function $f:\X  \to \power(\N^{<\omega})$
such that $C = f^{-1}[ {\WF}] $. 
\item (Assuming $\X = Y^\N$.) There is
 a set $R \subseteq Y^{<\omega}\times \N^{<\omega}$ such that 
$$ A  = \{ f\in Y^\N :  \ \exists g\in \N^\N \,\forall n\,\,
	(f\on n, g\on n) \in R \}$$
\end{enumerate}
\end{Fact}

The coanalytic sets are just the sets whose complement is analytic. 
Borel sets are of course both analytic and coanalytic, and the 
``Separation theorem'' states that the converse is true:

Let $A\subseteq \X$ be both analytic and coanalytic. 
Then $A$ is a Borel set.  

Analytic sets have the following closure properties:
\begin{Fact}\label{closure}
\begin{enumerate}
\item All Borel sets are analytic (and coanalytic).
\item The countable union or intersection of analytic sets is again
analytic.  Similarly, the countable union or intersection of
coanalytic sets is again analytic.
\item The continuous preimage of an analytic set is analytic. The
continuous preimage of a coanalytic set is coanalytic. 
\item The continuous image of an analytic set is analytic.  
(Note that the continuous image of a Borel set is in general not Borel.)
\item
In particular, if $C \subseteq \N^\N \times \N^\N$ is a Borel set, then the
set $\{f\in \N^\N: \exists g\in \N^\N\,
	 (f,g)\in C\}$ is analytic, and the set  
 $\{f\in \N^\N: \forall g\in \N^\N\, (f,g)\in C\}$ is coanalytic. 
\end{enumerate}
\end{Fact}
However, while the Borel sets are closed under complements, the
analytic sets are not.   There are coanalytic sets which are not
analytic, for example the set ${\WF}$. 

We call a set $D \subseteq Y^X$ ``complete coanalytic'' iff
\begin{enumerate}
\item $D$ is coanalytic
\item For any coanalytic set $C\subseteq Y^X$ there is a continuous
function 
$F:Y^X\to Y^X$ with $C= F^{-1}[D]$. 
\end{enumerate}

It is known that the set~$\WF$ 
is complete coanalytic.  In fact, $\WF$ is the
``typical'' coanalytic set: 

Let $D$ be coanalytic.  Then $D$ is complete coanalytic iff there is a 
continous function $F:\power(\N^{<\omega})\to Y^X$ with $\WF = F^{-1}[D]$.

Equivalently, $D$ is complete coanalytic iff there is a function as above
which is defined only on the set of trees.

The existence of coanalytic sets which are not analytic easily implies
that  a complete coanalytic set can never  be analytic. 

The following theorem should be read as ``analytic sets can never
reach $ \omega_1$.''

\begin{Fact}[Boundedness theorem]\label{bound}
\
\begin{enumerate}
\item  Every coanalytic set is the union of an increasing
$\omega_1$-chain of Borel sets. 
\item Let $\WF = \bigcup_{\alpha \in \omega_1} \WF_\alpha$ be an 
increasing union of Borel sets, and let $A \subseteq \WF $ be Borel
(or even analytic). \\
Then there is  $\alpha\in \omega_1 $ such that 
 $A\subseteq \WF_\alpha$. 
\end{enumerate}

\end{Fact}

\section{Clones below $\T_1$}
\label{t1}

\begin{Definition}
We fix a 1-1 function $p$ from $\N\times \N $ onto $ \N\setminus \{0\}$. 
Let $\chi_\Delta$ and  $\chi_\nabla $ be the characteristic functions
of $\Delta$ and $\nabla$, and let $p_\Delta:= p \cdot \chi_\Delta$,
i.e., $p_\Delta(x,y) = p(x,y)$ for $x>y$, and $=0$ otherwise. 

Similarly, let 
$p_\nabla:= p\cdot \chi_\nabla$.  
\end{Definition}

The following is clear: 
\begin{Fact} \
\begin{itemize}
\item $\chi_\nabla$ and $\chi_\Delta$ are canonical,
		 and  in $\T_1\cap \T_2$.
\item $p_\Delta $ and $p_\nabla$ are in $\T_1\setminus \T_2$, and 
are canonical. 
\item $p\notin \T_1\cup \T_2$.   In fact, the only 
clone containing $\oo 1 \cup \{p\}$ is $\fullclone$ itself. 
\end{itemize}
\end{Fact}

%
%
%
%
%

\begin{Theorem} $\T_1$ is generated by $\{p_\Delta\}\cup \oo 1$. 
\end{Theorem}

\begin{proof}
Let $\C$ be the binary clone generated by $\{p_\Delta\}\cup \oo 1$.
We will first find a function $q\in \C$ satisfying
\begin{enumerate}
\item $q$ is 1-1 on $\Delta$
\item $q(x,y) = Q(x)$ on $\nabla$,  for some 1-1 function $Q$
\item $q[\Delta]\cap q[\nabla] = \emptyset$. 
\end{enumerate}
Note that  any two  functions $q$, $q'$  satisfying these
properties will be equivalent, in the sense that there is a unary 
function $u$ with $q(x,y)=u(q'(x,y))$ for all $(x,y)\in \Delta\cup \nabla$.)


	\begin{figure}

	\setlength{\unitlength}{0.0004in}
	\begingroup\makeatletter\ifx\SetFigFont\undefined%
	\gdef\SetFigFont#1#2#3#4#5{%
		\reset@font\fontsize{#1}{#2pt}%
			\fontfamily{#3}\fontseries{#4}\fontshape{#5}%
			\selectfont}%
			\fi\endgroup%
{\renewcommand{\dashlinestretch}{30}
	\begin{picture}(5037,3264)(0,-10)
		\drawline(1575,3237)(1575,12)(5025,12)
		\drawline(1575,12)(4725,3162)
		\put(3750,1062){\makebox(0,0)[lb]{$p(x,y)$}}
	\put(2475,2412){\makebox(0,0)[lb]{$x$}}
	\end{picture}
}

	\begin{center} 
Properties of $q$ (simplified)
	\end{center}

	\end{figure}

	Let $P(x) = \max \{p(x,y): y\le x\}+1$, and let  
	$$  q(x,y):= p_\Delta(P(x), p_\Delta(x,y)). $$
	Note that this actually means $q(x,y) = 
	p(P(x), p_\Delta(x,y))$, as $P(x)>p_\Delta(x,y)$ for all $x,y$.
	So, $$
	q(x,y) = \KNUTHcases{ p(P(x), p(x,y)) & for  $x>y$\cr
		p(P(x),0)   & for  $x \le y$,\cr}$$

		So $  q$ satisfies (1)--(3), and $q\in \C$. 

		We now consider an arbitrary almost unary function  $f$, 
		say $f(x,y) < F(x)$ for all $x,y$.  Wlog we assume $f(x,y)>0$ for all
		$(x,y)$.

		Let $p':\N\times \N \to \N$ be a 1-1 function satisfying $p'(x,y)>x$ for
		all $x,y$.

Define
	$$
	\begin{aligned}
	f_1(x,y) &= \KNUTHcases{
		p'(F(x), f(x,y)) 
			& for $x> y$\cr
			F(x)     & for $x\le y$\cr
	}\\
f_2(x,y) &= \KNUTHcases{\rlap{0}
	\hphantom{	 p'(F(x), f(x,y)) }
	& for $x> y$\cr
		f(x,y)    & for $x\le y$\cr
}
\end{aligned}
$$
		\begin{figure}
		\setlength{\unitlength}{0.0004in}
		\begingroup\makeatletter\ifx\SetFigFont\undefined%
		\gdef\SetFigFont#1#2#3#4#5{%
			\reset@font\fontsize{#1}{#2pt}%
				\fontfamily{#3}\fontseries{#4}\fontshape{#5}%
				\selectfont}%
				\fi\endgroup%

{\renewcommand{\dashlinestretch}{30}
	\begin{picture}(9037,3264)(0,-10)
		\drawline(0075,3237)(0075,12)(3525,12)
		\drawline(0075,12)(3225,3162)
	\put(0975,2412){\makebox(0,0)[lb]{$F(x)$}}
	\put(1750,562){\makebox(0,0)[lb]{$f(x,y)$ (but $>x$)}}

	\drawline(5575,3237)(5575,12)(8925,12)
		\drawline(5575,12)(8725,3162)
	\put(6475,2412){\makebox(0,0)[lb]{$f(x,y)$}}
	\put(7250,562){\makebox(0,0)[lb]{$0$}}
	\end{picture}
}

\begin{center} Definitions of $f_1$ and $f_2$ (simplified)
	\end{center}

	\end{figure}

Then $f_1(x,y) = u_1(q(x,y))$ for some unary $u_1$, and 
$f_2(x,y) = u_2(p_\Delta(y+1, x)$ for some unary $u_2$.   
		So $f_1,f_2\in \C$.

		Let $f'(x,y):= p_\Delta(f_1(x,y),f_2(x,y))$. 
		Now $f_2(x,y)<   F(x) \le f_1(x,y)$ for all $x,y$, so
		$f'(x,y) = p(f_1(x,y),f_2(x,y))$.    

		As $f(x,y)$ can be recovered from
		the pair $(f_1(x,y), f_2(x,y))$, and hence also from $f'(x,y)$,  we 
		conclude that $f(x,y) = v(f'(x,y))$ for some unary $v$. Hence $f\in
		\C$. 
		\end{proof}

		\begin{Theorem} If $\C \subseteq \T_1$ is a binary clone containing $\oo 1$, 
		then either $\C= \T_1$, or $\C\subseteq \T_2$. 

		Hence:  $\T_1\cap \T_2$ is the unique coatom in the interval $[\oo 1,
		\T_1]$ of binary clones, and every binary 
		clone in this interval (except for $\T_1$ itself) is
		included in $ \T_1\cap \T_2$. 
		\end{Theorem}

\begin{proof} 
Assume $ \oo 1 \subseteq \C \subseteq \T_1$, but $\C\not \subseteq
\T_2$.  So let $f\in \C\setminus \T_2$.   So there are 1-1 unary
functions $u$ and $v$ such that $f(u(x), v(y))$ is canonical and 1-1
on $\Delta$ (or on $\nabla$). So wlog $f$ is canonical and 1-1 on $\Delta$.

Moreover, either $f$ is
symmetric, or $\ran(f\on \nabla)\cap \ran(f\on \Delta)= \emptyset$. 

In the first case, the function $f'(x,y):=f(2x, 2y+1)$ is 1-1 on all of
$\N\times \N$, so $\langle \{f'\}\cup \oo 1\rangle = \fullclone $, which contradicts our
assumption $\C\subseteq \T_1$. 

In the second case, we can find a
unary  function $u$ such that $$\forall x,y: \ u(f(2x,2y+1)) =
p_\Delta(x,y),$$ so $p_\Delta\in\C$, i.e., $C=\T_1$.

\end{proof}

Pinsker \cite{Pinsker:2004a} has analyzed the interval $(\T_1,\Pol(\T_1))$ of
(full) clones, and shown the following:   

\newcommand{\second}{\min^+}
\begin{Theorem}[Pinsker]  
Let $\second_n(x_1,\ldots, x_n):= x_{\pi (2)}$, where $\pi$ is any
permutation such that $x_{\pi(1)}\le x_{\pi(2)}\le \cdots \le
x_{\pi(n)}$. \\ (So $\second_2(x,y)=\max(x,y)$, and $\second_3(x,y,z)$
		is the median of $x,y,z$.) 

Then the clones $\M_n:=\langle \T_1\cup \{\second_n\}\rangle$ are all
distinct,
$$\T_1 \subseteq \cdots \subsetneq \M_5 
	\subsetneq \M_4 
\subsetneq \M_3 = \Pol(\T_2)
	\subsetneq \M_2 = \fullclone,   $$
	and every clone in the interval $[\T_1,\Pol(\T_1)]$ is equal to some
	$\M_n$. 
	\end{Theorem}

\begin{Remark} So $\M_4$ is a coatom in the 
 interval 
$[\oo 1 , \Pol(\T_1)]$ in the lattice of all clones.   It is also 
easy to see that $
\Pol(\T_1)\cap \Pol(\T_2) = \Pol(\T_1\cap \T_2) $ is another  coatom. 
\end{Remark}

	\begin{Fact} $\T_1$ is a Borel set. 
	\end{Fact}
	\begin{proof} The set $\T_1^{\rm x}:= \{f\in\oo 2:  \exists F \, \forall
		x,y:f(x,y)\le F(x)\}$ is apparently only $\Sigma_1^1$, but we can
		rewrite it as 
		\begin{equation*}
		\begin{split}
		\T_1^{\rm x} &= \{
			f\in\oo 2: \forall x\,\exists z\, \forall y: f(x,y)\le z \} \\
			&= 
			\bigcap_{x\in \N}   \, 
			\bigcup_{z\in \N} \, 
			\bigcap_{y\in \N}  \, 
			\bigcup_{t\le z} \, 
			\{f\in \oo 2: f(x,y) = t \},
			\end{split}
			\end{equation*}
			which is  $F_{\sigma\delta}$. 

			$\T_1^{\rm y}$ can be defined similarly, and $\T_1 = \T_1^{\rm x}\cup 
			\T_1^{\rm y}$. 

			\end{proof}

			\begin{Remark} Clearly, $\Pol(\T_1)$ is coanalytic. (See
\ref{closure}(5).) 
			By Pinsker's theorem, $\Pol(\T_1) = \langle \T_1\cup
			\{\second_3\}\rangle$  is finitely generated over
			$\T_1$, hence analytic and therefore even Borel.   An explicit
			Borel description can be found in \cite{Pinsker:2004a}.
			\end{Remark}

\section{Clones below $\T_2$}
\label{t2}

In the previous section we have seen: 
\begin{Theorem}
$\T_1 = \langle \oo 1 \cup \{p_\Delta\}\rangle$.
  Thus, $\T_1$ is finitely generated
over $\oo 1 $. 
\end{Theorem}

  The next theorem and its corollaries show that $\T_2$ is not finitely
generated over $\oo 1 $.

%
%
%
%
%
%
%

\begin{Fact} Let $B \subseteq \fullclone$ be a Borel or analytic set.   
Then $\langle B\rangle$ is analytic. 

Similarly, if $B \subseteq \oo 2 $ is a Borel or analytic set, then 
 $\langle B\rangle_{\oo 2}$ (the binary clone generated by $B$)
is analytic. 
\end{Fact}

\begin{Question}  Is there a Borel set $B$ (perhaps even a closed set?
      a countable set?  a set of the form $\oo 1 \cup \{f_1,\ldots, f_n\}$?) 
 such that $\langle B\rangle$ is
      not Borel? 
\end{Question}

\begin{Theorem}
 $\T_2$,  $\Pol(\T_2)$, $\T_1\cap \T_2$ and $\Pol(\T_1\cap \T_2)$
  are complete coanalytic sets.
\end{Theorem}

\begin{proof}
We will define a continuous map $F$ from the set of all trees $T
      \subseteq \N^{<\N}$ into $\T_1$ such that
\begin{quote}   for all $T $:  $T$ is wellfounded iff $F(T)\in \T_2$. 
\end{quote}

Let $\{s_n:n\in \N\}$ enumerate all finite sequences of natural
numbers, with $ s_k\vartriangleleft s_n \Rightarrow k<n$. 

  For any tree $T \subseteq \{s_n:n\in \N\}$ let $F(T)$ be
defined as follows: 
$$ F(T)(k,n) = \KNUTHcases { p(k,n) & if $k<n$ and $s_k,s_n\in T$,
      $s_k\vartriangleleft s_n$\cr
			0	& otherwise\cr
}
$$

Now if $A:=\{s_{n_1}, s_{n_2},\ldots \,\}$ is an infinite branch in $T$,
      then $F(T)\on\nabla_{A,A} $ is 1-1.  

Conversely:  Assume  $A = \{n_1<n_2< \cdots \}$, $B= \{m_1<m_2<\cdots \}$,
and $F(T)\on \nabla_{A,B}$ is 1-1.   \\
We claim that  $s_{n_1}\vartriangleleft s_{n_2}$.  Indeed, for any large
enough $k$ we have $F(T)(n_1,m_k)\not=0$, so 
 $s_{n_1}\vartriangleleft s_{m_k}$, and similarly 
 $s_{n_2}\vartriangleleft s_{m_k}$.
 So $s_{n_1}\vartriangleleft s_{n_2}$.  
\\
Similarly we get $s_{n_1}\vartriangleleft s_{n_2} \vartriangleleft 
s_{n_3}\vartriangleleft \cdots $. 
\end{proof}

\begin{Corollary}
$\Pol(\T_2)$, $\T_2$, $\T_2\cap \T_1$ are not countably generated over
      $\oo 1$. 
\end{Corollary}
\begin{proof} If $C$ is a countable set, then $C\cup \oo 1 $ is Borel,
so 
$\langle C\cup \oo 1 \rangle$
is analytic, hence not complete coanalytic.
\end{proof}

The well-known analysis of  coanalytic sets now gives the following: 
\begin{Theorem}
There is a sequence $(C_i:i\in \omega_1)$ of Borel clones such that:
      
$i<j$ implies $C_i\subsetneq  C_j$,  $\bigcup_{i\in \omega_1} C_i =
      \T_2$, and:\\ 

For every analytic clone $C\subseteq \T_2$ there is $i< \omega_1$ such
      that $C\subseteq C_i$. 

In other words: There is an increasing family of $\aleph_1$ many
analytic clones below $\T_2$ such that every analytic clone below
$\T_2$ is covered by a clone from the family. 

A similar representation can be found for $\Pol(\T_2)$, $\T_1\cap \T_2$, etc.

\end{Theorem}
\begin{proof}
By \ref{bound}(1), 
we can find an increasing family of Borel sets $(B_i:i<\omega_1)$ such
 that $\T_2=\bigcup_i B_i$. 
Clearly each clone $\langle B_i\rangle$ is analytic.  By the
boundedness theorem (\ref{bound}(2))
  we know that for all $i$ there is $j$ with
$\langle B_i\rangle \subseteq B_j$.   Let $h:\omega_1\to \omega_1$ be
continuous and strictly increasing such that $\forall i: \langle
B_i\rangle \subseteq B_{f(i)}$.     Now the family $\{B_i: f(i)=i\}$
is as desired. 
\end{proof}

\begin{Question} Find a nice cofinal family in $\{\C: \C \subsetneq
\T_2\}$.  I.e., a nice family $\F$ such that $\forall \C\subsetneq \T_2$
there is $\C'\in \F$ with $\C\subseteq \C'$.  

Since we already have a family that covers all analytic clone, this
question asks really:  which nonanalytic clones are there below $\T_2$?

\end{Question}

\begin{Question}
Can we get a family $B_i$ as in the theorem where each $B_i$ is
generated by a single function? 
\end{Question}

\begin{Question}
Analyze the interval $[ \T_2, \Pol(\T_2)]$. 
\end{Question}

%
%

\bibliographystyle{plain}
\bibliography{listb,listx,goldstrn,other}

\begin{thebibliography}{1}

\bibitem{Gavrilov:1965}
G.P. Gavrilov.
\newblock {Ueber funktionale Vollstaendigkeit in der abzaehlbar-wertigen
  Logik.}
\newblock {\em Probl. Kibernetiki}, 15:5--64, 1965.

\bibitem{GoSh:737}
Martin Goldstern and Saharon Shelah.
\newblock {Clones on regular cardinals}.
\newblock {\em Fundamenta Mathematicae}, 173:1--20, 2002.
\newblock {\tt arXiv:math.RA/0005273}

\bibitem{Kechris:1995}
Alexander~S. Kechris.
\newblock {\em Classical descriptive set theory}, volume 156 of {\em Graduate
  Texts in Mathematics}.
\newblock Springer-Verlag, New York, 1995.

\bibitem{Moschovakis:1980}
Yiannis~Nicholas Moschovakis.
\newblock {\em {Descriptive Set Theory}}, volume 100 of {\em Studies in Logic
  and the Foundations of Mathematics}.
\newblock North-Holland Publishing Company, Amsterdam New York Oxford, 1980.

\bibitem{Pinsker:2004a}
Michael Pinsker.
\newblock {Clones containing all almost unary functions}.
\newblock {\em Algebra Universalis}, to appear, 2004.
\newblock {\tt arXiv:math.RA/0401102}

\bibitem{PK:1979}
R.~P{\"o}schel and L.~A. Kalu{\v{z}}nin.
\newblock {\em Funktionen- und {R}elationenalgebren}, volume~15 of {\em
  Mathematische Monographien [Mathematical Monographs]}.
\newblock VEB Deutscher Verlag der Wissenschaften, Berlin, 1979.
\newblock Ein Kapitel der diskreten Mathematik. [A chapter in discrete
  mathematics].

\end{thebibliography}

\end{document}